\documentclass[leqno]{degm}

\newcommand{\Q}{{\mathbb{Q}}}
\newcommand{\F}{{\mathbb{F}}}
\newcommand{\C}{{\mathbb{C}}}
\newcommand{\Z}{{\mathbb{Z}}}

\newcommand{\Sym}{{\mathfrak{S}}}
\newcommand{\fB}{{\mathfrak{B}}}
\newcommand{\fg}{{\mathfrak{g}}}
\newcommand{\fh}{{\mathfrak{t}}}
\newcommand{\cE}{{\mathcal{E}}}
\newcommand{\cF}{{\mathcal{F}}}
\newcommand{\cG}{{\mathcal{G}}}
\newcommand{\cM}{{\mathcal{M}}}

\newcommand\Irr{\operatorname{Irr}}
\newcommand\Ind{\operatorname{Ind}}
\newcommand\Hom{\operatorname{Hom}}

\newtheorem{thm}{Theorem}[section]
\newtheorem{cor}[thm]{Corollary}
\newtheorem{prop}[thm]{Proposition}
\newtheorem{lem}[thm]{Lemma}

\theorembodyfont{\rmfamily}

\newtheorem{defn}[thm]{Definition}
\newtheorem{exmp}[thm]{Example}
\newtheorem{rem}[thm]{Remark}

\renewcommand{\leq}{\leqslant}
\renewcommand{\geq}{\geqslant}
\renewcommand{\atop}[2]{\genfrac{}{}{0pt}{}{#1}{#2}}

\contact[geck@desargues.univ-lyon1.fr]{Meinolf Geck, Institut Girard 
Desargues, Universit\'e Lyon 1, 21 av Claude Bernard, 69622 Villeurbanne 
cedex, France}

\begin{document}

\title{On the Schur indices of cuspidal unipotent characters}

\author{Meinolf Geck}

\maketitle

\begin{center} 
{\em To John Thompson on his $70$th birthday} 
\end{center}
\bigskip

\begin{abstract}
In previous work of Gow, Ohmori, Lusztig and the author, the Schur indices 
of all unipotent characters of finite groups of Lie type have been 
explicitly determined except for six cases in groups of type $F_4$, $E_7$ 
and $E_8$. In this paper, we show that the Schur indices of all cuspidal 
unipotent characters for type $F_4$ and $E_8$ are~$1$, assuming that the 
group is defined over a field of ``good'' characteristic. This settles
four out of the six open cases.  For type $E_7$, we show that the Schur 
indices are at most~$2$.
\end{abstract}

\begin{classification} Primary 20C15; Secondary 20G40
\end{classification}

\section{Introduction} \label{MGsec1}
By the work of Feit, Gow, Ohmori, Turull, $\ldots$, one expects that the
Schur index of an irreducible character of a finite simple group is $1$ 
or~$2$. This paper is a contribution to a solution of this problem, as far 
as the unipotent characters of a finite  group of Lie type are concerned.
By Lusztig \cite{LuRat}, the Schur indices of all rational-valued unipotent 
characters are generally~$1$. However, there are exceptions in the finite 
unitary groups (Ohmori \cite{Ohm1}) and the Ree groups ${^2\!F_4}(q)$ 
(\cite{myert03}) where the Schur indices are~$2$. In \cite{myert03}, we 
were also able to compute the Schur indices of some non-rational valued 
unipotent characters, but not all. This paper deals with some of  the 
remaining cases, which arise in  groups of type $F_4$, $E_7$ and $E_8$. 

Our principal tool are Kawanaka's generalized Gelfand--Graev characters 
\cite{Kaw2}. These are characters induced from certain unipotent subgroups. 
It is known that every unipotent character occurs with ``small'' 
multiplicity in such an induced character. We shall show here that 
certain generalized Gelfand--Graev representations can actually be realized
over $\Q$ (in type $F_4$ or $E_8$) or over a quadratic extension of $\Q$ 
(in type $E_7$). The method for doing this is inspired by Ohmori \cite{Ohm2},
who considered the case of ordinary Gelfand--Graev representations. 
Furthermore, using the explicit multiplicity formulas obtained by Kawanaka 
\cite{Kaw2} and Lusztig \cite{LuUS}, we check that the unipotent characters 
that we wish to consider all occur with multiplicity~$1$ in those generalized 
Gelfand--Graev representations. 

Combining the previous work by Ohmori, Lusztig and the author with the
results obtained in this paper, the known results on Schur indices of 
unipotent characters can be summarized as follows. Let $G$ be a simple 
algebraic group and $F\colon G \rightarrow G$ be an endomomorphism such 
that $G^F$ is a finite (twisted or untwisted) group of Lie type. Let 
$\rho$ be a unipotent  irreducible character of $G^F$. Then the Schur 
index of $\rho$ is given as follows. 
\begin{itemize}
\item[(1)] Assume that $\rho$ is cuspidal and rational-valued.
\begin{itemize}
\item[(a)] If $G^F$ is of type ${^2\!}A_{n-1}$ where $n=s(s+1)/2$ for 
some $s\geq 1$, then the Schur index of the unique cuspidal unipotent 
character of $G^F$ is $1$ if $n$ is congruent to $0$ or $1$ modulo $4$,
and $2$ otherwise; see Ohmori \cite{Ohm1}. 
\item[(b)] If $G^F$ is a Ree group of type ${^2\!}F_4$, then the Schur
index of the unique cuspidal unipotent character which occurs with even
multiplicity in all Deligne--Lusztig generalized characters $R_{T,1}$ 
is~$2$; see \cite{myert03}, Theorem~1.6. 
\item[(c)] In all other cases, the Schur index of $\rho$ is~$1$; see 
Lusztig \cite{LuRat}, Theorem~0.2. 
\end{itemize}
\item[(2)] Assume that $\rho$ is cuspidal and not rational-valued. Then the
Schur index of $\rho$ is~$1$, except possibly for the two cuspidal unipotent
characters in type $E_7$ (where the Schur index is at most $2$) or the two 
cuspidal unipotent characters with character field ${\Q}(\sqrt{-1})$ in 
type $E_8$ over a field of characterictic~$5$; see Table~1 and 
Corollary~1.5 in \cite{myert03}, and Sections~\ref{MGsec2}, \ref{MGsec3} 
in the present paper. See Gow \cite{Gow1}, p.~119, for groups of type 
${^2\!B}_2$ and ${^2\!}G_2$.
\item[(3)] Assume that $\rho$ is unipotent but not necessarily cuspidal. 
Then $\rho$ occurs with non-zero multiplicity in the Harish-Chandra induction
$R_L^G(\psi)$ where $L$ is an $F$-stable Levi complement in some $F$-stable
proper parabolic subgroup of $G$ and $\psi$ is a cuspidal unipotent character
of $L^F$. In this situation, the Schur index of $\rho$ equals that of 
$\psi$; see \cite{myert03}, Proposition~5.6.
\end{itemize}
Thus, there are two remaining cases to be dealt with:
\begin{itemize}
\item the characters $E_7[\pm \xi]$ in type $E_7$; 
\item the characters $E_8[\pm \sqrt{-1}]$ for a group
of type $E_8$ in characteristic~$5$.
\end{itemize}
For each of these cases, we show that the determination of the Schur 
index can be reduced to the problem of  computing explicitly the values 
of certain induced characters. The latter problems are discussed in 
\cite{mye7c} and \cite{hezard}.

\section{On the rationality of generalized Gelfand--Graev representations} 
\label{MGsec2}

The aim of this section is to show that certain generalized Gelfand--Graev
representations of a finite group of Lie type can be realized over $\Q$.
For this purpose, we have to recall in some detail the construction of
these representations. We shall freely use standard results and notations
concerning (connected reductive) linear algebraic groups and their Lie 
algebras (see Carter \cite{Ca2} and the references there).

Let $G$ be a connected reductive group defined over the finite field 
$\F_q$, with corresponding Frobenius map $F\colon G \rightarrow G$. We 
fix an $F$-stable Borel subgroup $B\subseteq G$ and write $B=U.T$ where 
$U$ is the unipotent radical of $B$ and $T$ is an $F$-stable maximal torus. 
Let $X=\Hom(T,k^\times)$ be the character group of $T$ and $\Phi\subseteq X$ 
be the root system of $G$ with respect to $T$. Then $B$ determines a 
positive system  $\Phi^+\subseteq \Phi$ and a corresponding set of simple 
roots $\Pi \subseteq \Phi^+$. We have 
\[ G=\langle T,X_\alpha\mid \alpha \in \Phi\rangle \qquad \mbox{and}
\qquad U=\prod_{\alpha \in \Phi^+} X_\alpha\]
where $X_\alpha$ is the root subgroup corresponding to $\alpha\in\Phi$.
(Here, it is understood that the product is taken over some fixed order
of the roots.) For each $\alpha\in \Phi$, let $\alpha^\vee\in 
Y=\Hom(k^\times,T)$ be the corresponding coroot. Given $\lambda\in X$ 
and $\mu \in Y$, let $\langle \lambda,\mu\rangle$ be the unique 
integer such that $\gamma^{\langle \lambda,\mu\rangle}=(\lambda\circ\mu)
(\gamma)$ for all $\gamma\in k^\times$. This defines a non-degenerate
pairing $\langle \;,\;\rangle \colon X\times Y\rightarrow\Z$; the matrix
\[ C:=\bigl(\langle \alpha,\beta^\vee\rangle\bigr)_{\alpha,\beta\in \Pi}\]
is the Cartan matrix of $G$. Now, for any $\alpha\in \Phi$, there is an 
isomorphism of algebraic groups $x_\alpha \colon k^+ \rightarrow X_\alpha$ 
such that 
\[ tx_\alpha(\xi)t^{-1}=x_\alpha(\alpha(t)\xi) \qquad \mbox{for
all $t\in T$ and $\xi \in k$}.\]
Since we will only need to consider this case, let us assume that $F$ is of 
split type, such that
\begin{alignat*}{2}
F(x_\alpha(\xi)) &=x_\alpha(\xi^q) \qquad &&\mbox{for all $\alpha\in \Phi$,
$\xi\in k$},\\
F(t) & =t^q \qquad &&\mbox{for all $t\in T$}.
\end{alignat*}
We consider the conjugacy classes of unipotent elements in $G$. For this
purpose, we assume throughout that the characteristic $p$ of $\F_q$ is a 
good prime for $G$. Recall that this means that $p$ is good for each 
simple factor involved in~$G$, and that the conditions for the various
simple types are as follows.
\[\begin{array}{rl} A_n: & \mbox{no condition}, \\
B_n, C_n, D_n: & p \neq 2, \\
G_2, F_4, E_6, E_7: &  p \neq 2,3, \\
E_8: & p \neq 2,3,5.  \end{array}\]
Then it is known (see Kawanaka \cite{Kaw2} and the references there;
see also Premet \cite{Premet}) that one can naturally attach to each 
unipotent class of $G$ a so-called {\em weighted Dynkin diagram}, that 
is, an additive map 
\[d:\Phi \rightarrow \Z \quad \mbox{such that}\quad 
d(\alpha)\in \{0,1,2\} \mbox{ for all $\alpha \in \Pi$}.\]
The assignment from unipotent classes to additive maps as above is 
injective, but not surjective in general. Complete lists of weighted
Dynkin diagrams for the various types of simple algebraic groups can
be found in Carter \cite{Ca2}, \S 13.1. Given such a weighted Dynkin diagram
$d$, the corresponding unipotent class is determined as follows. We set 
\[ L_d:=\langle T,X_\alpha \mid \alpha\in \Phi, d(\alpha)=0\rangle\qquad
\mbox{and}\qquad U_{d,i}:=\prod_{\atop{\alpha\in \Phi^+}{d(\alpha) \geq i}} 
X_\alpha\]
for $i=1,2,3,\ldots$. Then $P_d:=U_{d,1}.L_d$ is a parabolic subgroup of 
$G$, with unipotent radical $U_{d,1}$ and Levi complement $L_d$. By Kawanaka
\cite{Kaw2}, Theorem~2.1.1, there is a unique unipotent class $C$ in $G$ 
such that $C\cap U_{d,2}$ is dense in $U_{d,2}$; furthermore, $C\cap U_{d,2}$ 
is a single $P_d$-conjugacy class and we have
\[C_G(u)\subseteq P_d \qquad \mbox{for all $u\in C\cap U_{d,2}$}.\]
Then $C$ is the unipotent class attached to the given weighted Dynkin 
diagram~$d$.

In order to define the generalized Gelfand--Graev representation associated
with an element $u\in C\cap U_{d,2}^F$, we need to introduce some further 
notation. Let $\fg$ be the Lie algebra of $G$ over $k=\overline{\F}_q$. 
Then $\fg$ is also defined over $\F_q$ and we have a corresponding 
Frobenius map $F\colon \fg \rightarrow\fg$. Let $\fh\subseteq \fg$ be 
the Lie algebra of $T$. We have a Cartan decomposition
\[ \fg=\fh \oplus \bigoplus_{\alpha\in \Phi} ke_\alpha \qquad \mbox{where
$F(\fh)=\fh$ and $F(e_\alpha)=e_\alpha$ for all $\alpha \in \Phi$}.\]
We shall set $c_\alpha:=\kappa(e_{\alpha}^*,e_{\alpha})$ for any $\alpha
\in \Phi^+$, where $\kappa\colon \fg \times \fg \rightarrow k$ is a 
non-degenerate $G$-invariant, associative bilinear form and $x\mapsto x^*$ 
is an opposition $\F_q$-automorphism of $\fg$, that is, an automorphism 
such that $\fh^*=\fh$ and $e_\alpha^*\in {\F}_qe_{-\alpha}$ for all $\alpha
\in \Phi$ (see Kawanaka \cite{Kaw2}, \S 3.1). 

Finally, we fix a non-trivial 
homomorphism $\chi\colon \F_q^{+} \rightarrow \C^\times$. 

\begin{defn}[Kawanaka] \label{gggr} Consider a unipotent element $u\in C
\cap U_{d,2}^F$; write  
\begin{equation*}
u \in \Bigl(\prod_{\atop{\alpha \in \Phi^+}{d(\alpha)=2}} 
x_\alpha(\eta_\alpha)\Bigr)\cdot U_{d,3}^F \qquad \mbox{where 
$\eta_\alpha\in \F_q$}.\tag{$*$}
\end{equation*}
With this notation, we define a map 
$\varphi_u\colon U_{d,2}^F \rightarrow \C^\times$ by
\[ \varphi_u\Bigl(\prod_{\atop{\alpha\in \Phi^+}{d(\alpha) \geq 2}} x_\alpha
(\xi_\alpha)\Bigr)=\chi\Bigl(\sum_{\atop{\alpha\in \Phi^+}{d(\alpha)=2}}
c_\alpha\,\eta_{\alpha}\, \xi_{\alpha}\Bigr) \qquad \mbox{where 
$\xi_\alpha\in \F_q$}.\]
Since $\kappa(e_\alpha,e_\beta)=0$ unless $\alpha=-\beta$, this definition
coincides with Kawanaka's original definition in  \cite{Kaw2}, (3.1.5); see
also \cite{myphd}, Chapter~1. The map $\varphi_u$ actually is a group 
homomorphism, that is, a linear character of $U_{d,2}^F$. Inducing that 
character from $U_{d,2}^F$ to $G^F$, we obtain 
\[ \Ind_{U_{d,2}^F}^{G^F}\bigl(\varphi_u\bigr)=[U_{d,1}^F:
U_{d,2}^F]^{1/2}\cdot \Gamma_u\]
where $\Gamma_u$ is the character of the generalized Gelfand--Graev 
representation associated with $u$; see Kawanaka \cite{Kaw2}, (3.1.12). 
Note that the index $[U_{d,1}^F:U_{d,2}^F]$ is an even power of~$q$, so 
the square root exists.
\end{defn}

\begin{rem} \label{allorb} In the above setting, let us choose another
element $u'\in C\cap U_{d,2}^F$, 
\[u'\in \Bigl(\prod_{\atop{\alpha \in \Phi^+}{d(\alpha)=2}}
x_\alpha(\eta_\alpha')\Bigr)\cdot U_{d,3}^F \qquad \mbox{where 
$\eta_\alpha'\in \F_q$}.\]
Then we can also apply the above constructions to $u'$. Thus, we
obtain a corresponding generalized Gelfand--Graev character 
$\Gamma_{u'}$ such that 
\[ \Ind_{U_{d,2}^F}^{G^F}\bigl(\varphi_{u'}\bigr)=[U_{d,1}^F:
U_{d,2}^F]^{1/2}\cdot \Gamma_{u'}.\]
Using the defining formula for $\Gamma_u$, it is straightforward to
check that $\Gamma_u=\Gamma_{u'}$ if $u'$ and $u$ lie in the same 
$P_d^F$-orbit. (Indeed, if $u'=gug^{-1}$ where $g\in P_d^F$, then 
$\varphi_{u'}(x)=\varphi_u(g^{-1}xg)$ for all $x\in U_{d,2}^F$ and so the
induced characters are equal.) Now let $C$ be the $G$-conjugacy class
of~$u$. Then $C^F$ splits into orbits under the finite group $G^F$, and
it is well-known that these $G^F$-orbits are parametrized by the 
$F$-conjugacy classes of $C_G(u)/C_G(u)^\circ$. Since $C_G(u) \subseteq 
P_d$, a complete set of representatives for the $G^F$-orbits in $C^F$ can 
be found inside $C\cap U_{d,2}^F$. 
\end{rem}

Assume now that $G$ is simple of adjoint type and write $\Pi=\{\alpha_1,
\ldots,\alpha_l\}$ where $l$ is the rank of $G$. Then we have $X=
{\Z}\Pi$. Let $\{\omega_1\,\ldots,\omega_l\}$ be the dual basis of $Y$,
that is, we have $\langle \alpha_i,\omega_j\rangle=\delta_{ij}$ for
$1\leq i,j\leq l$. Now we can state the following result which generalizes 
the arguments in Ohmori \cite{Ohm2}, Lemma~2 and Proposition~1(i). (See 
also Example~\ref{ohm} below.) 

\begin{prop} \label{ohm1} Let $G$ be simple of adjoint type and assume 
that there exist integers $n_j\in \Z$ ($1\leq j \leq l$) such that
\[\sum_{j=1}^l n_j\, \langle\alpha,\omega_j\rangle=1 \quad \mbox{for 
all $\alpha\in \Phi^+$ such that $d(\alpha)=2$ and $\eta_{\alpha}\neq 0$},\]
where the coefficients $\eta_\alpha$ are defined as in Definition~\ref{gggr}
($*$). Then the character 
$[U_{d,1}^F:U_{d,2}^F]^{1/2}\cdot \Gamma_{u}$ can be realized over $\Q$. 
If, moreover, we also have 
\[\sum_{j=1}^l n_j\, \langle\alpha,\omega_j\rangle=0 \quad \mbox{for 
all $\alpha \in \Phi^+$ such that $d(\alpha)=0$},\] 
then $[U_{d,1}^F:U_{d,2}^F]^{1/2}\cdot \Gamma_{u'}$ can be realized over 
$\Q$, for every $u'\in C\cap U_{d,2}^F$. 
\end{prop}

\begin{proof} Let $\nu$ be a generator for the multiplicative group of 
$\F_q$. We claim that there exists an element $t\in T^F$ such that 
\begin{equation*}
\alpha(t)=\nu \quad \mbox{for all $\alpha\in \Phi^+$ with $d(\alpha)=2$
and $\eta_\alpha\neq 0$}.\tag{$\dagger$}
\end{equation*}
This is seen as follows. By \cite{Ca2}, Proposition~3.1.2, the map 
\[ h\colon \underbrace{k^\times \times \cdots\times  k^\times}_{\text{$l$ 
factors}} \rightarrow T, \qquad (x_1,\ldots,x_l)\mapsto\prod_{i=1}^l 
\omega_i(x_i),\] 
is an $\F_q$-isomorphism. In particular, we have $T^F=\{h(x_1,\ldots,x_l)
\mid x_i\in {\F}_q^\times\}$. 
Now set $t=h(\nu^{n_1},\ldots,\nu^{n_l})\in T^F$. Let $\alpha\in \Phi^+$ be 
such that $d(\alpha)=2$ and $\eta_\alpha\neq 0$. Then we have 
\[ \alpha(t)=\prod_{j=1}^l x_j^{\langle\alpha, \omega_j\rangle}=
\prod_{j=1}^l \nu^{n_j\langle\alpha, \omega_j\rangle}=
\nu^{\sum_{j=1}^l n_j\langle\alpha, \omega_j\rangle}=\nu,\]
as required. Thus, ($\dagger$) is proved.

Now let $H:=\langle t \rangle$. Then $H$ is a group of order $q-1$ which 
normalizes $U_{d,2}^F$. Let us induce $\varphi_u$ from $U_{d,2}^F$ to the
semidirect product $U_{d,2}^F.H$ and denote the induced character by 
$\psi_u$.  We will show that $\psi_u$ is a rational-valued irreducible 
character. First note that $\psi_u$ has non-zero values only on elements in 
$U_{d,2}^F$. Furthermore, by Mackey's formula, the restriction of $\psi_u$ 
to $U_{d,2}^F$ is given by $\sum_{h \in H} \varphi_u^h$ where $\varphi_u^h
(x)= \varphi_u (hxh^{-1})$ for all $x\in U_{d,2}^F$. Now let 
\[x:= \prod_{\atop{\alpha\in \Phi^+}{d(\alpha) \geq 2}} x_\alpha
(\xi_\alpha)\in U_{d,2}^F \qquad \mbox{for any $\xi_\alpha \in \F_q$}.\]
Then we have 
\[txt^{-1}=\prod_{\atop{\alpha\in \Phi^+}{d(\alpha) \geq 2}} tx_\alpha
(\xi_\alpha)t^{-1}=\prod_{\atop{\alpha\in \Phi^+}{d(\alpha) \geq 2}} 
x_\alpha (\alpha(t)\xi_\alpha)\]
and so, using ($\dagger$),
\[ \varphi_u^t(x) = \chi\Bigl(\sum_{\atop{\alpha\in 
\Phi^+}{d(\alpha)=2}} c_\alpha\,\eta_{\alpha}\, \alpha(t)
\xi_{\alpha}\Bigr)= \chi\Bigl(\nu\sum_{\atop{\alpha\in 
\Phi^+}{d(\alpha)=2}} c_\alpha\,\eta_{\alpha}\,\xi_{\alpha}\Bigr).\]
Similarly, for any $1\leq i\leq q-1$, we obtain 
\[ \varphi_u^{t^i}(x)=\chi(\nu^i \gamma_x) \qquad \mbox{where}
\qquad \gamma_x:=\sum_{\atop{\alpha\in \Phi^+}{d(\alpha)=2}} c_\alpha\,
\eta_{\alpha}\,\xi_{\alpha}\in {\F}_q.\]
This shows, first of all, that the characters $\varphi_u^{t^i}$ ($1\leq i
\leq q-1$) are pairwise different and, hence, the character $\psi_u$ is 
seen to be irreducible. Furthermore, we obtain  that 
\[ \psi_u(x)=\sum_{h \in H} \varphi_u^h(x)=\sum_{i=1}^{q-1} \varphi_u^{t^i}
(x)= \sum_{i=1}^{q-1} \chi(\nu^i\gamma_x)=\left\{\begin{array}{cl} q-1 &
\mbox{ if $\gamma_x=0$},\\ -1 & \mbox{ if $\gamma_x\neq 0$}.\end{array}
\right.\]
In particular, this shows that $\psi_u(x)\in \Z$ for all $x\in U_{d,2}^F$.
Thus, $\psi_u$ is a rational-valued irreducible character of $U_{d,2}^F.H$.
In order to show that $\psi_u$ can be realized over $\Q$, we note that
(again by Mackey's formula), the restriction of $\psi_u$ to $H$ is the
character of the regular representation of $H$. Hence the trivial character
$1_H$ occurs with multiplicity~$1$ in that restriction. By Frobenius
reciprocity, $\psi_u$ itself occurs with multiplicity~$1$ in the character
of $U_{d,2}^F.H$ obtained by inducing $1_H$. Since the latter character 
clearly is realized over $\Q$, a standard argument on Schur indices
shows that $\psi_u$ can be realized over $\Q$ (see Isaacs \cite{Isa},
Corollary~10.2). Finally, by the transitivity of induction, the character 
\[ \Ind_{U_{d,2}^F}^{G^F}\bigl(\varphi_u\bigr)=[U_{d,1}^F:U_{d,2}^F]^{1/2}
\cdot \Gamma_u\]
can be realized over $\Q$, as claimed.

Now assume that the integers $n_j$ satisfy the additional equations
described above and consider an arbitrary element $u'\in C\cap U_{d,2}^F$. 
Thus, we have $u'=gug^{-1}$ where $g\in P_d$. Let 
\[u'\in \Bigl(\prod_{\atop{\alpha \in \Phi^+}{d(\alpha)=2}}
x_\alpha(\eta_\alpha')\Bigr)\cdot U_{d,3}^F \qquad \mbox{where 
$\eta_\alpha'\in \F_q$}.\]
We claim that if $\eta_\alpha'\neq 0$, then $\alpha=\alpha_1+\beta$
where $d(\beta)=0$, $d(\alpha_1)=2$ and $\eta_{\alpha_1}\neq 0$. To see 
this, it is enough to consider the special case where $g\in T$ or where
$g\in X_\beta$ with $d(\beta)\geq 0$. The assertion is clear if $g\in T$.
On the other hand, if $g\in X_\beta$ where $d(\beta)\geq 0$, then 
Chevalley's commutator relations show that $gug^{-1}$ will lie in a 
product of root subgroups $X_{i\beta+\alpha_1}$ where $d(\alpha_1)\geq 2$
and $i=0,1,2,\ldots$, and the assertion also follows. We can now conclude 
that the integers $n_j$ automatically have the property that 
\[\sum_{j=1}^l n_j\, \langle\alpha,\omega_j\rangle=1 \quad 
\mbox{for all $\alpha\in \Phi^+$ such that $d(\alpha)=2$ and 
$\eta_{\alpha}'\neq 0$}.\]
Thus, we can apply the above arguments to $u'$ as well and conclude that
the character $[U_{d,1}^F:U_{d,2}^F]^{1/2}\cdot \Gamma_{u'}$ can be realized
over $\Q$.
\end{proof}

\begin{exmp} \label{ohm} (Cf.\ Ohmori \cite{Ohm2}, p.~151.) Assume that
$G$ is simple of adjoint type. The weighted Dynkin diagram corresponding to 
the conjugacy class of regular unipotent elements in $G$ is the map $d_0
\colon \Phi \rightarrow \Z$ such that $d_0(\alpha_i)=2$ for all $1\leq i 
\leq l$. In this case, we have 
\[P_{d_0}=B,\qquad L_{d_0}=T\qquad\mbox{and}\qquad U_{d_0,1}=U_{d_0,2}=U;\]
furthermore, $u=\prod_{i=1}^l x_{\alpha_i}(\eta_i)$, where $\eta_i\in 
\F_q^\times$, is a representative in that class. Now the system of equations
in Proposition~\ref{ohm1} is simply given by 
\[ \sum_{j=1}^l n_j\langle \alpha_i ,\omega_j\rangle=1 \qquad \mbox{for all
$1\leq i \leq l$}.\]
Thus, we have $n_j=1$ for $1\leq j \leq l$. Consequently, $\Gamma_u$ 
(which is just an ordinary Gelfand--Graev character) can be realized 
over  $\Q$.
\end{exmp}

\begin{exmp} \label{dist} Assume that $G$ is simple of adjoint type and 
that $d$ is a weighted Dynkin diagram of the following special type: There 
exists some $i_0\in \{1,\ldots, l\}$ such that 
\[d(\alpha_{i_0})=2 \quad \mbox{and}\quad d(\alpha_i)=0 \quad 
\mbox{for $i\in \{1,\ldots,l\}\setminus \{i_0\}$}.\]
We claim that, in this case, the generalized Gelfand--Graev character
$\Gamma_u$ can be realized over $\Q$, for any $u \in C\cap U_{d,2}^F$.
This is seen as follows. First note that $U_{d,1}=U_{d,2}$. Now let 
$\eta_\alpha$ be defined as in Definition~\ref{gggr}($*$) and consider the 
two sets of equations in Proposition~\ref{ohm1}. The second set of equations 
reads 
\[ \sum_{j=1}^l n_j \langle \alpha_i,\omega_j\rangle=0 \quad 
\mbox{for all $i \in \{1,\ldots,l\}\setminus \{i_0\}$}.\]
This yields $n_i=0$ for all $i\neq i_0$. Then the first set of equations
reads 
\[ n_{i_0} \langle \alpha,\omega_{i_0}\rangle =1
\quad \mbox{for all $\alpha \in \Phi^+$ such that $d(\alpha)=2$ and 
$\eta_\alpha\neq 0$}.\]
However, writing any such $\alpha$ as a linear combination of simple
roots, we see that $\alpha_{i_0}$ always occurs with coefficient~$1$ 
and so the above equations hold with $n_{i_0}=1$. Thus, the systems
of equations in Proposition~\ref{ohm1} have a unique solution.
\end{exmp}

We close this section with a general remark concerning the Schur
indices of unipotent characters. This remark shows that it will be enough 
to consider the unipotent characters of simple groups of adjoint type.

\begin{rem} \label{center} Assume that $G$ is simple of adjoint type.
Let $G_1$ be an algebraic group over ${\F}_q$ such that $G_1/Z(G_1)$ is
simple of the same type as $G$. We denote the Frobenius map on $G_1$ again 
by $F$. Then we have a surjective homomorphism of algebraic groups $\varphi 
\colon G_1\rightarrow G$ which is defined over ${\F}_q$ and such that
$\ker(\varphi)=Z(G_1)$. Let $\chi$ be a unipotent character of $G^F$. Denote 
by $\chi_1$ the character of $G_1^F$ which is obtained by first restricting 
$\chi$ to $\varphi(G_1^F)$ and then pulling back to $G_1^F$ via $\varphi$. 
By \cite{DeLu}, Proposition~7.10, the map $\chi\mapsto \chi_1$ defines
a bijection between the unipotent characters of $G^F$ and $G_1^F$, 
respectively. One can express this by saying that the unipotent characters 
are ``insensitive'' to the center of $G$.  

For a unipotent character
$\chi\in \Irr(G^F)$, denote by ${\Q}(\chi)={\Q}(\chi(g)\mid g \in G^F)$ the 
character field of $\chi$ and by $m_{\Q}(\chi)$ the Schur index of $\chi$; 
we use similar notations for a unipotent character of $G_1^F$. We claim that
\[{\Q}(\chi)={\Q}(\chi_1)\qquad\mbox{and}\qquad m_{\Q}(\chi)=m_{\Q}(\chi_1).\]
\end{rem}

\begin{proof} First note that we clearly have ${\Q}(\chi_1)\subseteq 
{\Q}(\chi)$. To show that we have equality, let $E\supseteq {\Q}$ be a 
finite Galois extension such that any unipotent character of $G^F$ and 
of $G_1^F$ can be realized over $E$. Given any $\tau \in \mbox{Gal}(E/\Q)$ 
and any unipotent character $\chi \in \Irr(G^F)$, we denote by $\chi^\tau$ 
the irreducible character obtained by algebraic conjugation; it is still 
a unipotent character of $G^F$ (see the remarks in \cite{myert03}, (5.1)). 
Now the map $\chi \mapsto \chi_1$ certainly is compatible with field 
automorphisms. Thus, we have $\chi_1^\tau=(\chi^\tau)_1$. Since the map 
$\chi\mapsto \chi_1$ is a bijection, we conclude that $\chi_1^\tau=\chi_1$ 
if and only if $\chi^\tau= \chi$. This implies that ${\Q}(\chi)=
{\Q}(\chi_1)$, as required.

Now consider the Schur indices. First note that we certainly have 
$m_{\Q}(\chi_1)\leq m_{\Q}(\chi)$. (Indeed, if $\chi$ can be realized
over some extension field of $\Q$, then $\chi_1$ is automatically realized 
over the same field.) To show the reverse inequality, we argue as follows.
The character $\chi_1$ can be realized over an extension field $E_1
\supseteq \Q$ such that $[E_1:\Q]=m_{\Q}(\chi_1)$. We now regard $\chi_1$ 
as a character of $\varphi(G^F)$. By Frobenius reciprocity, $\chi$ occurs
with multiplicity~$1$ in the character obtained by inducing $\chi_1$ from 
$\varphi(G^F)$ to $G^F$. Since that induced character can also be realized 
over $E_1$, a standard argument on Schur indices (\cite{Isa}, 
Corollary~10.2) shows that $m_{\Q} (\chi)\leq m_{\Q} (\chi_1)$, as required.
\end{proof}

\section{Simple groups of type $G_2$, $F_4$ and $E_8$} \label{MGsec3}
Throughout this section, let $G$ be simple of (adjoint) type $G_2$, 
$F_4$ or $E_8$, defined over the finite field ${\F}_q$ with corresponding 
Frobenius map $F \colon G \rightarrow G$. Then $F$ is of split type, in 
the sense of the previous section. The cuspidal unipotent characters of 
$G^F$ and the previously known results on their Schur indices are listed 
in Table~\ref{tab1}.

\begin{table}[htbp] \caption{Schur indices of cuspidal unipotent 
characters in type $G_2$, $F_4$ and $E_8$ (notation of Carter \cite{Ca2},
\S 13.9; see also \cite{myert03}, Table~1)} \label{tab1}
\begin{center}
$\renewcommand{\arraystretch}{1.2}\begin{array}{ccl}\hline 
\mbox{Type}&\mbox{Character} & \mbox{Schur index}\\\hline 
G_2  & G_2[1], G_2[-1] & 1 \quad\mbox{(see Lusztig \cite{LuRat})}\\
     & G_2[\theta], G_2[\theta^2] & 1 \quad \mbox{(see 
     \cite{myert03})}\\\hline
F_4       & F_4^I[1], F_4^{II}[1], F_4[-1]& 1 \quad\mbox{(see 
Lusztig \cite{LuRat})} \\ & F_4[i],F_4[-i] & ? 
\quad\mbox{($1$ in char. $2$, $3$; see \cite{myert03})}\\
           & F_4[\theta],F_4[\theta^2] & 1 \quad 
	   \mbox{(see \cite{myert03})}\\ \hline
E_8 & E_8^I[1],E_8^{II}[1],E_8[-1]& 1 
\quad\mbox{(see Lusztig \cite{LuRat})}\\
& E_8[-\theta],E_8[-\theta^2], & ?\quad\mbox{($1$ in char. 
$2$, $3$, $5$; see \cite{myert03})}\\ & E_8[\theta],E_8[\theta^2] & ? 
\quad\mbox{($1$ in char. $2$, $3$, $5$; see 
\cite{myert03})}\\ & E_8[i],E_8[-i] & ? \quad\mbox{($1$ in char.
$2$, $3$; see \cite{myert03})}\\ & E_8[\zeta],E_8[\zeta^2], E_8[\zeta^3],
E_8[\zeta^4]& 1 \quad \mbox{(see \cite{myert03})}\\ \hline 
\multicolumn{3}{c}{\mbox{$i:=\sqrt{-1}$, $\quad\theta:=\exp(2\pi i/3)$,
$\quad\zeta:=\exp(2\pi i /5)$.}}\end{array}$
\end{center}
\end{table}

Now let us assume that $G$ is defined over a field of good characteristic.
All the cuspidal unipotent characters of $G^F$ have the same ``unipotent 
support'' in the sense of Lusztig \cite{LuUS}; see also \cite{gema}, 
Proposition~4.2. Hence there is a unique unipotent class $C_0$ in $G$ such
that all the cuspidal unipotent characters of $G^F$ have unipotent support 
$C_0$. Let $e_0=\dim \fB_u$ (the dimension of the variety of Borel 
subgroups containing an element $u\in C_0$). Then, by \cite{gema}, 
Theorem~3.7, $q^{e_0}$ is the exact power of $q$ dividing the degreess of 
the unipotent characters of $G^F$ having unipotent support $C_0$. By 
inspection of the tables in \cite{Ca2}, Chapter~9, we find that $e_0=1$,
$4$ or $16$ for $G$ of type $G_2$, $F_4$ or $E_8$, respectively. Furthermore,
for $G$ of type $G_2$ or $F_4$, the class $C_0$ is uniquely determined
by this condition. In type $E_8$, we use the additional information
that $C_G(u)/C_G(u)^\circ$ must be non-trivial (see once more the formula 
in \cite{gema}, Theorem~3.7). Then $C_0$ is uniquely determined. The weighted 
Dynkin diagram $d_0$ associated to $C_0$ is specified in Table~\ref{tabwd}. 
Furthermore, we have $C_G(u)/C_G(u)^\circ \cong \Sym_3$, $\Sym_4$ or
$\Sym_5$ for $G$ of type $G_2$, $F_4$ or $E_8$, respectively.

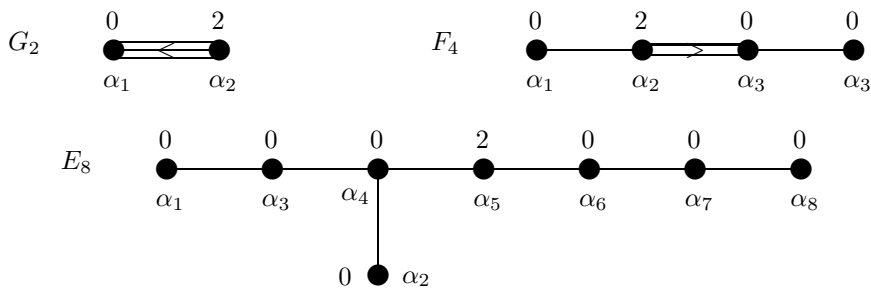
\begin{table}[htbp] \caption{Weighted Dynkin diagrams for the unipotent
supports of cuspidal unipotent characters in type $G_2$, $F_4$ and $E_8$}
\label{tabwd}
\begin{center} 
\begin{picture}(330,105)
\put(  0,90){$G_2$}
\put( 40,90){\circle*{8}}
\put( 40, 87){\line(1,0){40}}
\put( 40, 90){\line(1,0){40}}
\put( 40, 93){\line(1,0){40}}
\put( 56, 87.5){$<$}
\put( 80, 90){\circle*{8}}
\put( 36, 76){$\alpha_1$}
\put( 37, 98){$0$}
\put( 76, 76){$\alpha_2$}
\put( 77, 98){$2$}
\put(160, 90){$F_4$}
\put(200, 90){\circle*{8}}
\put(200, 90){\line(1,0){40}}
\put(240, 90){\circle*{8}}
\put(240, 88){\line(1,0){40}}
\put(240, 92){\line(1,0){40}}
\put(256, 87.5){$>$}
\put(280, 90){\circle*{8}}
\put(280, 90){\line(1,0){40}}
\put(320, 90){\circle*{8}}
\put(196, 76){$\alpha_1$}
\put(197, 98){$0$}
\put(236, 76){$\alpha_2$}
\put(237, 98){$2$}
\put(276, 76){$\alpha_3$}
\put(277, 98){$0$}
\put(316, 76){$\alpha_3$}
\put(317, 98){$0$}
\put( 20, 45){$E_8$}
\put( 60, 45){\circle*{8}}
\put( 60, 45){\line(1,0){40}}
\put(100, 45){\circle*{8}}
\put(100, 45){\line(1,0){40}}
\put(140, 45){\circle*{8}}
\put(140, 45){\line(1,0){40}}
\put(180, 45){\circle*{8}}
\put(180, 45){\line(1,0){40}}
\put(220, 45){\circle*{8}}
\put(220, 45){\line(1,0){40}}
\put(260, 45){\circle*{8}}
\put(260, 45){\line(1,0){40}}
\put(300, 45){\circle*{8}}
\put(140, 45){\line(0,-1){40}}
\put(140,  5){\circle*{8}}
\put( 56, 31){$\alpha_1$}
\put( 57, 53){$0$}
\put( 96, 31){$\alpha_3$}
\put( 97, 53){$0$}
\put(126, 33){$\alpha_4$}
\put(137, 53){$0$}
\put(176, 31){$\alpha_5$}
\put(177, 53){$2$}
\put(216, 31){$\alpha_6$}
\put(217, 53){$0$}
\put(256, 31){$\alpha_7$}
\put(257, 53){$0$}
\put(296, 31){$\alpha_8$}
\put(297, 53){$0$}
\put(149,  2){$\alpha_2$}
\put(125,  1){$0$}
\end{picture}
\end{center}
\end{table}

Now let us fix a cuspidal unipotent character $\rho$ of $G^F$. By Lusztig, 
\cite{LuUS}, Theorem~11.2, there exists some $u\in C_0\cap U_{d_0,2}^F$
such that $\rho$ occurs with ``small'' multiplicity in $\Gamma_u$. In
the present situation, Kawanaka actually obtained an explicit multiplicity
formula. To state that formula, we need to recall some facts about the 
parametrization of unipotent characters; see Lusztig \cite{Lu4}, Main 
Theorem~4.23. 

Let $\cF_0$ be the unique family of unipotent characters which contains 
the cuspidal unipotent characters. Let $\cG_0$ be the finite group attached 
to that family; we have $\cG_0\cong C_G(u)/C_G(u)^\circ$ in the present
situation. Let $\cM_0$ be the set of all pairs $(x,\sigma)$ where $x\in 
\cG_0$ (up to conjugacy) and $\sigma\in \Irr(C_{\cG_0}(x))$. Then we have 
a bijection 
\[\cM_0\rightarrow\cF_0,\qquad (x,\sigma) \mapsto \rho_{(x,\sigma)}.\]
(which satisfies some further conditions which we do not need to recall
here). On the other hand, there is a well-defined ``split'' element $u_0
\in C_0\cap U_{d_0,2}^F$; see Shoji \cite{Sh1}, Remark~5.1. Having fixed 
that element, there is a canonical parametrization of the $G^F$-conjugacy
classes contained in $C_0^F$ by the conjugacy classes of $\cG_0$. We shall 
denote by $\Gamma_{u_y}$ the generalized Gelfand--Graev character with 
respect to a unipotent element $u_y\in C_0\cap U_{d_0,2}^F$ in the 
$G^F$-class parametrized by $y\in \cG_0$. Having fixed this notation, we 
can now state the following result.

\begin{thm}[Kawanaka] \label{mult} Recall that $G$ is a group of type 
$G_2$, $F_4$ or $E_8$ in good characteristic. Let $\langle\rho_{(x,\sigma)},
\Gamma_{u_y}\rangle$ be the multiplicity of $\rho_{(x,\sigma)}$ in the
generalized Gelfand--Graev character $\Gamma_{u_y}$, where $x,y\in \cG_0$
and $\sigma \in \Irr(C_{\cG_0}(x))$. Then we have 
\[\langle\rho_{(x,\sigma)},\Gamma_{u_y}\rangle=\left\{\begin{array}{cl}
\sigma(1)&\quad\mbox{if $x=y \pmod{\cG_0}$},\\0&\quad
\mbox{otherwise}.  \end{array}\right.\]
\end{thm}

\begin{proof} Kawanaka obtained an explicit formula for the values of 
all the generalized Gelfand--Graev characters of $G^F$; see \cite{Kaw2}, 
Corollary~3.2.7 and Lemma~3.3.10. Furthermore, in \cite{Kaw2}, \S 4.2,
he obtained explicit formulas for the values of the unipotent characters
on unipotent elements. (In the latter reference, it is actually assumed that
the characteristic is large. But this hypothesis is only used in reference
to results on the Green functions of $G^F$, which are now known to hold
in general; see Shoji \cite{Sh2}, Theorem~5.5.) The above multiplicity 
formula is obtained by a simple formal rewriting of Kawanaka's formulas. 
In large characteristic, that formula could also be deduced from Lusztig's 
more general results in \cite{LuUS}. 
\end{proof} 

\begin{cor} \label{MGfe} Recall that $G$ is a group of type $G_2$, $F_4$ 
or $E_8$ in good characteristic. Then every cuspidal 
unipotent character of $G^F$ has Schur index~$1$.
\end{cor}

\begin{proof} By inspection of the tables in \cite{Ca2}, \S 13.9, all 
cuspidal unipotent characters of $G^F$ are labelled by pairs 
$(x,\sigma)\in \cM_0$ where $\sigma(1)=1$. Thus, by Theorem~\ref{mult}, 
every such character occurs with multiplicity~$1$ in some generalized 
Gelfand--Graev character attached to an element in $C_0\cap U_{d_0,2}^F$. 
Now we notice that the weighted Dynkin diagrams in 
Table~\ref{tabwd} are of the type considered in Example~\ref{dist}.
Thus, all the generalized Gelfand--Graev characters associated with
unipotent elements in $C_0\cap U_{d_0,2}^F$ can be realized over $\Q$.  
Hence, the assertion follows by a standard argument on Schur indices
(see Isaacs \cite{Isa},  Corollary~10.2).
\end{proof} 

The above result also covers some of the cases (at least in good 
characteristic) already dealt with by Lusztig \cite{LuRat} or the 
author \cite{myert03}. In combination with the previously known results
listed in Table~\ref{tab1}, we see that there is only one case left that 
remains to be considered:
\begin{center} 
{\em the characters $E_8[\pm i]$ for a group of type $E_8$ in 
characteristic~$5$.}
\end{center}
We shall now outline a strategy for dealing with this case. That strategy 
is analogous to that in \cite{myert03}, (6.5). For the remainder of this 
section, let $G$ be a simple group of type $E_8$ over $\F_q$ where $q$ 
is odd and $p\equiv 1 \bmod 4$. Let us label the simple roots of $G$ as in 
Table~\ref{tabwd} and let $\{\omega_1,\ldots, \omega_8\}$ be the dual basis 
of $Y$. Then we have $T=\{h(x_1,\ldots,x_8) \mid x_i\in k^\times\}$ as in 
the proof of Proposition~\ref{ohm1}. Now let $\nu$ be a generator of the 
multiplicative group of $\F_q$ and consider the semisimple element 
\[ s:=h(1,1,1,1,1,\nu^{(q-1)/4},1,1) \in T^F.\]
Thus, $s$ is an element of order $4$. Its centralizer is given by 
$G_1=\langle T,X_\alpha\mid \alpha \in \Phi_1\rangle$ where 
$\Phi_1=\{\alpha \in \Phi\mid \alpha(s)=1\}$. It is readily checked that a
system of simple roots of $\Phi_1$ is given by $\Pi_1=(\Pi\setminus 
\{\alpha_6\})\cup \{-\alpha_0\}$, where 
\[\alpha_0=2\alpha_1+3\alpha_2+ 4\alpha_3 +6\alpha_4+5\alpha_5 +4\alpha_6+
3\alpha_7+2\alpha_8\]
is the unique highest root in $\Phi$. The root system $\Phi_1$ is of
type $D_5\times A_3$. Now, by Lusztig \cite{LuCS}, Proposition~21.3, there 
is a well-defined unipotent class $C_1\subseteq G_1$ which supports a 
``cuspidal local system''. That unipotent class is determined explicitly
by \cite{LuSp}, Corollary~4.9, as far as the $D_5$-factor is concerned.
(We necessarily have the class of regular unipotent elements in the
$A_3$-factor; see \cite{LuCS}, \S 18.) Using \cite{Ca2}, \S 13.9, we see 
that the weighted Dynkin diagram $d_1$ of $C_1$ is given as follows:

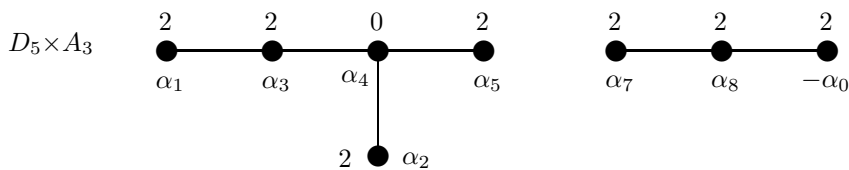
\begin{center}
\begin{picture}(330,60)
\put(  0, 45){$D_5{\times} A_3$}
\put( 60, 45){\circle*{8}}
\put( 60, 45){\line(1,0){40}}
\put(100, 45){\circle*{8}}
\put(100, 45){\line(1,0){40}}
\put(140, 45){\circle*{8}} 
\put(140, 45){\line(1,0){40}}
\put(180, 45){\circle*{8}}
\put(230, 45){\circle*{8}}
\put(230, 45){\line(1,0){40}}
\put(270, 45){\circle*{8}}
\put(270, 45){\line(1,0){40}}
\put(310, 45){\circle*{8}}
\put(140, 45){\line(0,-1){40}}
\put(140,  5){\circle*{8}}
\put( 56, 31){$\alpha_1$}
\put( 57, 53){$2$}
\put( 96, 31){$\alpha_3$}
\put( 97, 53){$2$}
\put(126, 33){$\alpha_4$}
\put(137, 53){$0$}
\put(176, 31){$\alpha_5$}
\put(177, 53){$2$}
\put(226, 31){$\alpha_7$}
\put(227, 53){$2$}
\put(266, 31){$\alpha_8$}
\put(267, 53){$2$}
\put(300, 31){$-\alpha_0$}
\put(307, 53){$2$}
\put(149,  2){$\alpha_2$}
\put(125,  1){$2$}
\end{picture}
\end{center}
By \cite{Ca1}, Theorem~11.3.2, we can explicitly identify an adjoint group 
of type $D_5$ with a $10$-dimensional orthogonal group. Determining the
Jordan normal form of the matrix corresponding to the element 
\[ u_1:=x_{\alpha_1}(1)x_{\alpha_5}(1)x_{\alpha_2}(1)
x_{\alpha_3+\alpha_4}(1)x_{\alpha_4+\alpha_5}(1)\cdot x_{\alpha_7}(1)
x_{\alpha_8}(1)x_{-\alpha_0}(1),\]
we see that $u_1\in C_1\cap U_{d_1,2}^F$.
Now, we have $C_{G_1} (u_1)/C_{G_1}^\circ(u_1)\cong {\Z}/4{\Z}$; see 
\cite{myert03}, Table~2, and the references there. Let $\psi$ be one 
of the two faithful irreducible characters of this group and define 
$\Gamma_{(C_1,\psi)}$ as in \cite{LuUS}, (7.5). Thus, $\Gamma_{(C_1,\psi)}$ 
is a certain linear combination (where the coefficients involve the values 
of $\psi$) of the generalized Gelfand--Graev characters of $G_1^F$ 
associated with the various unipotent elements in $C_1\cap U_{d_1,2}^F$.
Now we can state the following result.  

\begin{prop} \label{e8char5} In the above setting, assume that the following
conditions hold:
\begin{itemize}
\item[(a)] The values of $\Gamma_{(C_1,\psi)}$ are given by Lusztig's 
formula \cite{LuUS}, Corollary~7.6.  (This is known to be true if the 
characteristic is ``large''.) 
\item[(b)] The assertion in \cite{LuRem}, Theorem~0.8, holds for cuspidal
character sheaves in $G$. (This is known to be true if the characteristic
is good; see \cite{Sh2}, Theorem~4.1.)
\end{itemize}
Then the Schur index of the cuspidal unipotent character $E_8[\pm i]$ is~$1$.
\end{prop}

\begin{proof}  We begin by showing that the generalized Gelfand--Graev
characters of $G_1^F$ which are associated with the unipotent elements in 
$C_1\cap U_{d_1,2}^F$ can be realized over $\Q$. For this purpose, we use
an argument analogous to that in Proposition~\ref{ohm1}.

We shall consider a system of equations as in Proposition~\ref{ohm1} (but 
with respect to the unipotent element $u_1\in G_1$ specified above):
\[\sum_{j=1}^8 n_j\, \langle\alpha,\omega_j\rangle=\left\{\begin{array}{cl}
1 & \;\mbox{for $\alpha\in \{\alpha_1,\alpha_5,\alpha_2, \alpha_3+\alpha_4,
\alpha_4+\alpha_5,\alpha_7,\alpha_8,-\alpha_0\}$},\\
0 & \; \mbox{for all $\alpha \in \Phi^+$ such that $d_1(\alpha)=0$}. 
\end{array}\right.\]
One readily checks that this system has a unique solution:
\[ (n_1,n_2, \ldots,n_8)= (1,1,1,0,1,-5, 1, 1).\]
Thus, $t:=h(\nu,\nu, \nu, 1,\nu,\nu^{-5}, \nu,\nu)\in T^F$ is a 
semisimple element of order $q-1$. 

Now let $u\in C_1\cap U_{d_1,2}^F$. Then, as in the proof of 
Proposition~\ref{ohm1}, we see that $\alpha(t)=\nu$ (if $\alpha$ is 
involved in the expresssion of $u$ as a poduct of root subgroup elements) 
or $\alpha(t)=1$ (if $d_1(\alpha)=0$). Let $H=\langle t\rangle$ and 
denote by $\psi_u$ the character obtained by inducing $\varphi_u$ from
$U_{d_1,2}^F$ to $U_{d_1,2}^F.H$. Following exactly the same kind of 
arguments as in the proof of Proposition~\ref{ohm1}, we see that $\psi_u$
is an irreducible character which can actually be realized over $\Q$.
Consequently, the generalized Gelfand--Graev character $\Gamma_u$ of 
$G_1^F$ can also be realized over $\Q$, as claimed. 

We now follow the line of argument in \cite{myert03}, (6.5).  First note
that the character values of $E_8[\pm i]$ generate the field ${\Q}(i)$;
see \cite{myert03}, Table~1. Thus, we must show that $E_8[\pm i]$ can be
realized over ${\Q}(i)$. Let $\Sigma$ be the $G$-conjugacy class of $su_1$. 
Then $\Sigma$ supports a ``cuspidal local system'' $\cE$; see 
\cite{myert03}, Table~2 and the references there. Using \cite{Lu4}, Main 
Theorem~4.23, the explicit description of the Fourier matrix for type $E_8$, 
and our assumption (b), we see that 
\[ \langle E_8[\pm i],\chi_{(\Sigma,\cE)}\rangle=\frac{1}{4}\varepsilon
\qquad \mbox{where $\varepsilon\in \C$, $|\varepsilon|=1$};\]
here, $\chi_{(\Sigma,\cE)}$ is a characteristic function of $(\Sigma,\cE)$
and $\langle\;,\;\rangle$ denotes the usual hermitian product on the space of
complex-valued class functions on $G^F$. On the other 
hand, since $C_{G_1}(u_1)/C_{G_1}^\circ(u_1)\cong {\Z}/4{\Z}$, the set $C_1
\cap U_{d_1,2}^F$ splits into four classes in $G_1^F$, with representatives 
$u_1,u_2, u_3,u_4\in C_1\cap U_{d_1,2}^F$ say. Then, our assumption (a) 
implies that we have  
\[ 4\,\chi_{(\Sigma,\cE)}=\sum_{r=1}^4 \sum_{\lambda \in \Irr(Z_1)} 
\xi_{r,\lambda}\, \Ind_{Z_1\times U_{d_1,2}^F}^{G^F} \bigl(\lambda 
\boxtimes \varphi_{u_r}\bigr)\]
where $Z_1:=Z(G_1)=\langle s\rangle$ and $\xi_{r,\lambda} \in {\Z}[i]$; 
see the analogous equation in \cite{myert03}, (6.5), and note that 
$Z(G_1)^\circ=\{1\}$. Now, since $Z_1\cap H=\{1\}$, it is easily shown 
that 
\[ \Ind_{Z_1\times U_{d_1,2}^F}^{(Z_1\times U_{d_1,2}^F).H}
\bigl(\lambda \boxtimes \varphi_{u_r}\bigr)= \lambda \boxtimes 
\Ind_{U_{d_1,2}^F}^{U_{d_1,2}^F.H}\bigl(\varphi_{u_r}\bigr)=\lambda
\boxtimes \psi_{u_r}\]
can be realized over ${\Q}(i)$. Thus, we have $4\,\chi_{(\Sigma,\cE)}=
\sum_\alpha r_\alpha \rho_\alpha$ where $r_\alpha \in {\Z}[i]$ and 
$\rho_\alpha$ are characters of $G^F$ which can be realized over 
${\Q}(i)$. Since $\langle E_8[\pm i], 4\chi_{(\Sigma,\cE)}\rangle= 
\varepsilon$ has absolute value~$1$, we conclude that the greatest common 
divisor of the multiplicities $\langle E_8[\pm i],\rho_\alpha\rangle$ must 
be~$1$. Hence, a standard argument on Schur indices (\cite{Isa}, 
Corollary~10.2) implies that $E_8[\pm i]$ can be realized over 
${\Q}(i)$, as desired.
\end{proof}

In principal, the hypothesis on $\Gamma_{(C_1,\psi)}$ can be verified by
an explicit computation for a group of type $D_5 \times A_3$ (in the same 
way as Kawanaka \cite{Kaw2} verified that hypothesis for groups of 
exceptional type). The details will be discussed in \cite{hezard}.

%
\section{Simple groups of type $E_7$} \label{MGsec4}

Let $G$ be a simple group of adjoint type $E_7$ in good characteristic. 
There are precisely two cuspidal unipotent characters of $G^F$ denoted by 
$E_7[\pm \xi]$ where $\xi=\sqrt{-q}$; see
the table in \cite{Ca2}, \S 13.9. By \cite{myert03}, Table~1, the
character values of $E_7[\pm \xi]$ lie in the field ${\Q}(\xi)$ and the 
Schur index of $E_7[\pm \xi]$ is $1$, $2$ or $4$. Furthermore, by 
\cite{myert03}, Example~6.4, we already know that the Schur index is
$1$ if $p\not\equiv 1 \bmod 4$ or if $q$ is not a square, where $p$ is 
the characteristic of $\F_q$. Thus, the remaining task is to determine 
the Schur index when $q$ is a square and $p\equiv 1 \bmod 4$.

Let $C_0$ be the common unipotent support of the two cuspidal unipotent
characters. Arguing as in the previous section (that is, checking the 
exact power of $q$ dividing the degree of $E_7[\pm \xi]$ and comparing 
with the tables in \cite{Ca2}, Chapter~9), we find that the corresponding
weighted Dynkin diagram $d_0$ is given as follows. 

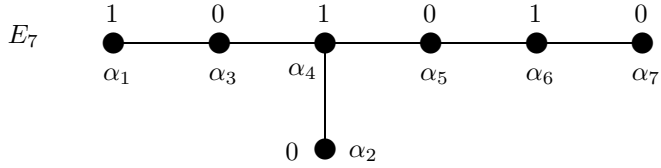
\begin{center}
\begin{picture}(280,60)
\put( 20, 45){$E_7$}
\put( 60, 45){\circle*{8}}
\put( 60, 45){\line(1,0){40}}
\put(100, 45){\circle*{8}}
\put(100, 45){\line(1,0){40}}
\put(140, 45){\circle*{8}} 
\put(140, 45){\line(1,0){40}}
\put(180, 45){\circle*{8}}
\put(180, 45){\line(1,0){40}}
\put(220, 45){\circle*{8}}
\put(220, 45){\line(1,0){40}}
\put(260, 45){\circle*{8}}
\put(140, 45){\line(0,-1){40}}
\put(140,  5){\circle*{8}}
\put( 56, 31){$\alpha_1$}
\put( 57, 53){$1$}
\put( 96, 31){$\alpha_3$}
\put( 97, 53){$0$}
\put(126, 33){$\alpha_4$}
\put(137, 53){$1$}
\put(176, 31){$\alpha_5$}
\put(177, 53){$0$}
\put(216, 31){$\alpha_6$}
\put(217, 53){$1$}
\put(256, 31){$\alpha_7$}
\put(257, 53){$0$}
\put(149,  2){$\alpha_2$}
\put(125,  1){$0$}
\end{picture}
\end{center}

\begin{lem} \label{lem2} Assume that $q$ is an even power of $p$. Let 
$u\in C_0\cap U_{d_0,2}^F$ and recall the definition of the linear
character $\varphi_u$ from Section~\ref{MGsec2}. Then there exists a 
field $K\supseteq \Q$ such that $[K:\Q]\leq 2$ and 
$\Ind_{U_{d_0,2}^F}^{G^F} \bigl(\varphi_{u} \bigr)$ can be realized
over $K$.
\end{lem}

\begin{proof} The following argument is inspired by Ohmori \cite{Ohm2},
Proposition~1. 

By Mizuno \cite{Miz}, Lemma~28, a representative $u_0\in C_0
\cap U_{d_0,2}^F$ is given by
\[u_0= x_{20}(1)x_{21}(1)x_{24}(1)x_{28}(1)x_{30}(1),\]
where the subscripts correspond to the following roots in $\Phi^+$:
\begin{alignat*}{2}
20&: \alpha_1+\alpha_2+\alpha_3+\alpha_4,&\qquad 21&: 
\alpha_1+\alpha_3+\alpha_4+\alpha_5,\\
24&: \alpha_2+\alpha_4+\alpha_5+\alpha_6,&\qquad 28&: 
\alpha_2+\alpha_3+2\alpha_4+\alpha_5,\\30&: \alpha_3+\alpha_4+
\alpha_5+\alpha_6+\alpha_7. &\qquad&
\end{alignat*}
We consider a set of equations similar to that in Proposition~\ref{ohm1}:
\[\sum_{j=1}^7 n_j\, \langle\alpha,\omega_j\rangle=\left\{\begin{array}{cl}
2 & \;\mbox{for $\alpha\in \Phi^+$ labelled by $20$, $21$, $24$, $28$ 
and $30$}, \\ 0 & \; \mbox{for all $\alpha \in \Phi^+$ such that
$d_0(\alpha)=0$}. \end{array}\right.\]
One readily checks that this system has a unique solution:
\[ (n_1,n_2,\ldots,n_7)=(1,0,0,1,0,1,0).\]
Now let $\nu$ be a generator of the multiplicative group of $\F_p\subseteq
\F_q$. Since $q$ is an even power of $p$, we can find a square root of 
$\nu$ in $\F_q$, which we denote by $\nu^{1/2}\in \F_q$. Then the element
\[ t:=h(\nu^{n_1/2},\ldots,\nu^{n_7/2})=h(\nu^{1/2},1,1,\nu^{1/2},1,
\nu^{1/2},1)\in T^F\]
has order $2(p-1)$ and, as in the proof of Proposition~\ref{ohm1}, we see
that $\alpha(t)=\nu$ for all $\alpha\in \Phi^+$ such that $d_0(\alpha)=2$
and $\eta_\alpha\neq 0$ where the coefficients $\eta_\alpha$ are defined 
as in Definition~\ref{gggr}($*$) (with respect to the unipotent element
$u\in C_0\cap U_{d_0,2}^F$). Now let $H=\langle t\rangle$. Then $H$ is a 
group of order $2(p-1)$ which normalizes $U_{d_0,2}^F$. Let us induce 
$\varphi_u$ from $U_{d_0,2}^F$ to the semidirect product $U_{d_0,2}^F.H$ 
and denote the induced character by $\psi_u$. Arguing as in the proof of 
Ohmori \cite{Ohm2}, Proposition~1(ii), we see that $\psi_u$ is the sum of 
two rational-valued irreducible characters $\psi_1$ and $\psi_2$, where
one of them, $\psi_1$ say, can be realized over $\Q$. By the Brauer--Speiser 
theorem (see \cite{CR2}, 74.27), the character $\psi_2$ has Schur index at 
most~$2$. Thus, there exists a field $K \supseteq \Q$ such that $[K:\Q]
\leq 2$ and $\psi=\psi_1+\psi_2$ can be realized over $K$. Consequently, 
by the transitivity of induction, the character
\[ \Ind_{U_{d_0,2}^F}^{G^F}\bigl(\varphi_u\bigr)=
\Ind_{U_{d_0,2}^F.H}^{G^F}\bigl(\psi_u\bigr)\]
can also be realized over $K$.
\end{proof}

\begin{rem} \label{rem3} We have remarked at the beginning of this section
that it is enough to compute the Schur index of $E_7[\pm \xi]$ in the case
where $q$ is a square and $p\equiv 1\bmod 4$. In this case, the character
values of $E_7[\pm \xi]$ lie in the field ${\Q}(\sqrt{-1})$. 

Now Ohmori \cite{Ohm2}, p.~154, points out that the character $\psi_2$ 
arising in the above proof has non-trivial local Schur indices at $\infty$ 
and at the prime $p$. This implies that $K$ is not contained in ${\Q}
(\sqrt{-1})$ if $p\equiv 1\bmod 4$. Indeed, if we had $K\subseteq {\Q}
(\sqrt{-1})$, then $\psi_2$ could be realized over ${\Q}(\sqrt{-1})$. 
Since $p\equiv 1 \bmod 4$, the field ${\Q}(\sqrt{-1})$ is contained in 
the field of $p$-adic numbers. Hence, $\psi_2$ could be realized over 
the latter field, contradicting the fact that the local 
Schur index at~$p$ is not one.
\end{rem}

\begin{cor} \label{Me7} Recall that $G$ is a group of type $E_7$ in
good characteristic. Then the two cuspidal unipotent characters of 
$G^F$ have Schur index at most~$2$.
\end{cor}

\begin{proof} We have already remarked at the beginning of this section
that the Schur index of $E_7[\pm \xi]$ is $1$ if $q$ is not a square. So
let us now assume that $q$ is an even power of $p$. Let $u_0\in C_0\cap 
U_{d_0,2}^F$. Then it is known that $C_G(u_0)/C_G(u_0)^\circ \cong \Z/2\Z$ 
(see the tables in \cite{Ca2}, \S 13.1). Thus, $C_0^F$ splits into two 
classes in the finite group $G^F$.  Let $u_1\in C_0\cap U_{d_0,2}^F$ be a 
representative of the second $G^F$-conjugacy class contained in $C_0^F$. 
Combining \cite{gema}, Theorem~3.7 and Remark~3.8, we have that
\[ \langle E_7[\pm \xi],\Gamma_{u_0}\rangle +\langle E_7[\pm \xi],
\Gamma_{u_1}\rangle=1.\]
(This could also be deduced from Kawanaka's explicit multiplicity
formulas \cite{Kaw2}.) Thus, $E_7[\pm \xi]$ occurs with multiplicity~$1$
in $\Gamma_{u}$ for some $u\in C_0\cap U_{d_0,2}^F$. Consequently, we
have
\[ \Bigl\langle E_7[\pm \xi], \Ind_{U_{d_0,2}^F}^{G^F}\bigl(\varphi_{u}
\Bigr)\Bigr\rangle=[U_{d_0,1}^F: U_{d_0,2}^F]^{1/2}=:m_0, \]
where the number $m_0$  on the right hand side is odd since we are in good
characteristic. Let $K\supseteq \Q$ be as in Lemma~\ref{lem2}. Then a 
standard result on Schur indices (see Isaacs \cite{Isa}, Lemma~10.4) shows 
that the Schur index of $E_7[\pm \xi]$ divides $[K:\Q]\, m_0$. However, as 
we have already remarked at the beginning of this section, the Schur index 
of $E_7[\pm \xi]$ must be $1$, $2$ or $4$. Thus, since $m_0$ is odd, we can
now conclude that the Schur index of $E_7[\pm \xi]$ divides $[K:\Q] \leq 2$.
\end{proof}

\begin{rem} \label{final}
The Schur index of $E_7[\pm \xi]$ equals $1$ or $2$, according to
whether  
\[ \Bigl\langle E_7[\pm \xi],\Ind_{U_{d_0,2}^F.H}^{G^F}\bigl(\psi_1\bigr)
\Bigr\rangle \mbox{ is odd or even},\]
where $\psi_1$ is defined in the proof of Lemma~\ref{lem2}. Thus, the problem
is reduced to the computation of the above scalar product, and this will be
discussed in \cite{mye7c}.
\end{rem}

\bigskip\noindent
{\bf Acknowledgements.} I wish to thank the organizers of the ``Finite
Groups 2003'' conference, C.~Y.~Ho, P.~Sin, P.~H.~Tiep and A.~Turull, for 
the invitation and the University of Gainesville for its hospitality.

\makelastpage
\end{document}